\documentclass[12pt]{article}
\usepackage[margin=.75in]{geometry}                
\usepackage{amsfonts}
\usepackage{amsmath}
\usepackage{amsthm}
\usepackage{amssymb}
\usepackage{tikz}
\usepackage{fancyhdr}
\usepackage{enumerate}
\usepackage{breakurl}

\usepackage{mathrsfs}  
\usepackage{colonequals}
\usepackage[hyphens,spaces,obeyspaces]{url}
\usepackage{hyperref}
\usepackage[boxed, section]{algorithm}
\usepackage{algpseudocode}
\usepackage[color = lime]{todonotes}
\usepackage{graphicx}
\usepackage[numbers]{natbib}

\theoremstyle{definition}
\theoremstyle{definition}\newtheorem{eg}{Example}
\theoremstyle{plain}
\theoremstyle{definition}
\theoremstyle{definition}
\theoremstyle{definition}
\theoremstyle{definition}
\theoremstyle{definition}
\theoremstyle{definition}
\theoremstyle{plain}\newtheorem{thm}{Theorem}
\theoremstyle{definition}
\newcommand{\comment}[1]{}

\theoremstyle{plain}
\theoremstyle{definition}
\theoremstyle{definition}
\theoremstyle{definition}
\theoremstyle{plain}

\newcommand{\C}{\mathbb{C}}

\newcommand{\Ff}{\mathbb{F}}

\newcommand{\Z}{\mathbb{Z}}

\newcommand{\Q}{\mathbb{Q}}

\newcommand{\red}{\operatorname{red}}

\newcommand{\p}{^\prime}

\renewcommand{\d}[1]{\ensuremath{\operatorname{d}\!{#1}}}

\title{Computing rational points on rank 0 genus 3 hyperelliptic curves}

\author{Mar\'ia In\'es de Frutos-Fern\'andez and Sachi Hashimoto}

\date{}

\begin{document}

\maketitle

\begin{abstract}
 	
 We compute rational points on genus $3$ odd degree hyperelliptic curves $C$ over $\mathbb{Q}$ that have Jacobians of Mordell-Weil rank $0$. The computation applies the Chabauty-Coleman method to find the zero set of a certain system of $p$-adic integrals, which is known to be finite and include the set of rational points $C(\Q)$. We implemented an algorithm in Sage to carry out the Chabauty-Coleman method on a database of $5870$ curves.
 	
\end{abstract}

\section{Introduction}

Given a curve $C$ of genus $g \geq 2$ defined over $\mathbb{Q}$, Mordell's conjecture, proved by Faltings \cite{Faltings} implies that the set of rational points $C(\mathbb{Q})$ is finite. Our goal is to explicitly compute this finite set of points in the case of a genus $3$ hyperelliptic curve with Jacobian having Mordell-Weil rank $0$. For a general curve, computing the set of its rational points is an unsolved problem. However, there are techniques available when the Mordell-Weil rank $r$ of the Jacobian is smaller than the genus. In this case it is often possible, for a specific curve, to provably find its rational points.

The rational point-finding techniques applied here for curves of small Mordell-Weil rank date back to an idea of Chabauty \cite{Chabauty}, who proved in 1941 that when the rank $r$ is less than the genus $g$, fixing a prime of good reduction $p$, the set $C(\mathbb{Q}_p) \cap \overline{J(\mathbb{Q})}$ is finite, and hence $C(\Q)$ is, where $\overline{J(\mathbb{Q})}$ denotes the closure in the $p$-adic topology of $J(\mathbb{Q})$. In the 1980s, Coleman made Chabauty's idea effective: he gave an upper bound on the number of rational points \cite{EffectiveChab}, now known as the Chabauty--Coleman bound. Coleman's upper bound comes from counting the zeros of $p$-adic (Coleman) integrals on $J(\mathbb{Q}_p)$ that vanish on $\overline{J(\mathbb{Q})}$. This bound was refined by Stoll \cite{Stoll}. We exhibit some curves that reach the Stoll bound. In 2010, Balakrishnan, Bradshaw, and Kedlaya \cite{BBK} gave a practical algorithm to compute these $p$-adic Coleman integrals on odd degree hyperelliptic curves and implemented the algorithm in Sage \cite{sagemath}.

Recently, Sutherland \cite{g3database} constructed a database of genus 3
hyperelliptic curves of discriminant up to $10^7$ using the approach of \cite{g2curves}. Armed with a decent upper bound on the number of rational points on a curve, a practical algorithm for computing these Coleman integrals that vanish on the rational points, and a large number of curves to experiment on, one might ask: how effective is Coleman's method for determining rational points?

For example, Balakrishnan, Bianchi, Cantoral-Farf\'an, \c{C}iperiani, and Etropolski \cite{CCexp} carried out computations to find rational points on 16,977 genus 3 odd degree hyperelliptic curves with rank 1 Jacobian.  They found that in most cases, the Chabauty-Coleman method picks up exactly the rational points, or the rational points and the Weierstrass points that are defined over $\Q_p$. In the remaining cases, they were able to explain the existence of the extra points either by torsion, linearity of the integral, or by extra endomorphisms of the Jacobian.

In order to satisfy the Chabauty-Coleman hypothesis that $r< g$ we need $ r = 0, 1,$ or $2$, but in the case of $r = g-1$, one will likely find many extra $p$-adic but non-rational points when applying the Chabauty-Coleman method, and in general the Chabauty-Coleman set will be larger than the set of rational points. More work would be needed to determine the set of rational points; for example, one could implement a Mordell-Weil sieve. For these reasons, we focused on the case of $g = 3$ and $ r = 0$. In this case, we computed the rational points in Magma \cite{magma}, up to height $10^5$ and verified these points using Chabauty-Coleman calculations. Most of the time, in $3070$ curves, the only extra points in the Chabauty-Coleman set were Weierstrass points that were defined over $\Q_p \setminus \Q$. In $17$ cases, we picked up points on the curve giving rise to torsion points of higher order on the Jacobian. These points are defined over quadratic fields where the prime $p$ splits. 

In Section 2, we present background on Coleman integration: we define the Coleman integral and list its properties, and we also discuss explicit Coleman integration. In Section 3, we discuss our implementation of the Chabauty--Coleman method for computing rational points on hyperelliptic curves of genus $3$ with rank $0$ Jacobian. Finally, in Section 4, we present our results and give an overview of the rational points and extra $\Q_p$-points found on the database of curves. We discuss several interesting examples in detail, including a curve with an $18$-torsion point on its Jacobian and a curve that is sharp for Stoll's bound.

\subsection*{Acknowledgments}

We would like to thank Jennifer Balakrishnan for suggesting these computations, and for many helpful discussions throughout. We would also like to thank Alex Best for several productive conversations about Coleman integration. The second author was supported by a Clare Boothe Luce grant while working on this project.

\section{Background on Coleman integration}

In this section we will define the basics of Coleman integration and necessary background for the rest of the paper. For a more thorough introduction to the rigid analytic geometry and theory behind Coleman integrals, we refer the reader to \cite{Torsion}, whose exposition we have loosely followed herein. Other good background on the Chabauty-Coleman method can be found in \cite{McCallumPoonen} and \cite{Wetherell}.

Fix $K$ a field that is complete with respect to a non-archimedean absolute value and fix some embedding of $K$ into $\C_p$. Let $R$ be the ring of integers of $K$.

Let $C/ R$ be a smooth connected affinoid curve over $K$ and set $C_K = C \times_R K$. Let $F$ be the residue field of $K$, and $\widetilde{C}_K$ be $C \times_R F$. There is a reduction map $\red:C_K \to \widetilde{C}_K$. We say that $C$ is \textit{has good reduction over $K$} if $\widetilde{C}_K$ is smooth. We call the preimage of a point under the reduction map $\red^{-1}(x)$ for $x \in \widetilde{C}_K$ the \emph{residue disc of $x$}. These residue discs partition $C_K(K)$.

Let $C$ have good reduction, $\omega$ be a holomorphic one-form on $C$, and $P, Q \in C(K)$ be $K$-rational points. Coleman defined an integral $\int_P^Q \omega$ in the following way. Fix a Frobenius lift $\phi: C_K \to C_K$, that is, a morphism of rigid analytic varieties which reduces to (relative) Frobenius on $\widetilde{C}_K$.

\begin{thm}[{\cite[Theorem~2.1]{Torsion}}]
\label{ColemansTheorem1}

	Suppose there is a polynomial $\mathcal{P}(T) \in \C_p[T]$ which does not vanish on any root of unity, and that \[ \mathcal{P}(\phi^*) \omega \text{ is exact.}\] Then there is a function $f_\omega$ on $C(\C_p)$ which is analytic on each residue disc such that $df_\omega = \omega$ and $ \mathcal{P}(\phi^*)(f_\omega)$ is analytic. The function $f_\omega$ is unique up to a constant and is independent of $\phi$ and $\mathcal{P}$.
	
\end{thm}

In this way, Coleman provides a way to compute a locally analytic antiderivative of $\omega$. Coleman proves, in several corollaries, properties of the antiderivative, which we collect here in one theorem:

\begin{thm}[{\cite[Theorem~2.3, Proposition~2.4]{Torsion}}]

	The integral satisfies the following properties:
	
	\begin{enumerate} 
		\item Linearity: $\int_P^Q (\alpha \omega + \beta \omega^ {\prime}) = \alpha \int_P^Q \omega + \beta \int_P^ Q \omega^\prime$.
		\item  Additivity: $\int_P^R \omega = \int_P^Q \omega + \int_Q^R\omega.$
		\item Change of variables: if $C^{\prime}$ is a curve and $\phi: C \to C^{\prime}$ a rigid analytic map between wide opens then $\int_P^Q \phi^* \omega = \int_{\phi(P)}^{\phi(Q)} \omega.$ For example, $\phi$ can be taken to be a lift of the $p$th power Frobenius.
		\item  Fundamental theorem of calculus: $\int_P^Q df = f(Q) - f(P)$.

		\item For any divisor $D\colonequals \sum_i Q_i - P_i$ of degree zero on $C$, the integral $\int_D \omega = \sum_i\int_{P_i}^{Q_i} \omega$ is well-defined, and $\int_D \omega = 0$ when $D$ is principal.

		\item Galois equivariance: if $\sigma$ is an automorphism of $\C_p$, then $\left( \int_P^Q \omega \right)^\sigma = \int_{\sigma(P)} ^{\sigma(Q)} \omega^\sigma$. 

	\end{enumerate}
	
\end{thm}

Then Coleman shows that we can use this integration theory to obtain arithmetic information about points.

Suppose now that $C_K$ is the base change of a curve $C$ defined over $\Q$ 
and that the Jacobian $J$ of $C/\Q$ has Mordell-Weil rank $r < g$. 
Fix $\iota: C \to J$ to be an embedding given by choosing a $\Q$-rational basepoint $b \in C(\Q)$ and mapping $P$ to $[P -b]$. Denote by $J_K$ the base change of $J$ to $K$.
 The embedding $\iota$ induces an isomorphism of vector spaces $\iota^*:H^0(J_K,\Omega^1) \to H^0(C_K, \Omega^1)$ which is independent of the chosen basepoint $b$ (Proposition 2.2 of \cite{MilneJacobian}).
 
We have functionals \[\lambda_\omega(D) = \int_{D} \omega, \text{ for } \omega \in H^0 (J_K, \Omega^1). \]
If $D = \sum_i [Q_i - P_i]$, then $\lambda_\omega$ can also be thought of as functionals on $C_K$ by
$\lambda_\omega (D) = \sum_i \int_{P_i}^{Q_i} \iota^*\omega.$  

Note that $H^0 (J_K, \Omega^1)$ is $g$-dimensional but the Jacobian $J(\Q)$ has rank $r < g$. Fix a basis $\omega_0, \dots, \omega_{g-1}$ for the space of holomorphic differentials on $J_K$. 
Define functionals \[\lambda_i (D) = \int_{D} \omega_i : J(\Q) \to K.\]
Consider all linear relations between the restrictions of the $\lambda_i$ to $J(\Q)$:
there are at least $g-r$ independent linear relations and they are of the form 
\begin{equation*}\label{lin-rel}
 \sum_{i=0}^{g-1} \alpha_i \int_{D}  \omega_i =0, \text{ with } \alpha_i \in K.
\end{equation*}

Let $S$ be the set of functionals $\mu(D) := \sum_{i=0}^{g-1} \alpha_i \int_{D}  \omega_i$ from  $J_K$ to $K$ such that $\mu(D) = 0$ for all $D \in J(\Q)$.
Via the map $\iota$, we can define $\mu(z) := \sum_{i=0}^{g-1} \alpha_i \int_{\iota(z)} \iota^* \omega_i$ for any $z \in C(K)$. Note that in particular $\mu(z) = 0$ for every $z$ in $C(\Q)$.

By Theorem \ref{ColemansTheorem1}, the functionals $\mu$ in $S$ are locally analytic, and we can use bounds on the numbers of zeros of $p$-adic power series to get control over the number of zeros for $\mu$ on each residue disc. In particular, there are at most finitely many common zeros of any non-zero linear combination of functionals. Hence the set of rational points $C(\Q)$ is contained in the finite set of points satisfying $\mu(z)=0$ for every $\mu$ in $S$.

Coleman's proof suggests an effective method for computing $C(\Q)$: first, compute the linear relations that the functionals satisfy. Then, compute the finite zero set of the functionals on each residue disc. This second step requires an algorithm for computing the Coleman integral to some $p$-adic precision.

Balakrishnan, Bradshaw, and Kedlaya give a practical algorithm for computing $\int_P^Q \omega$, on an odd-degree hyperelliptic curve $C$, where $P$ and $Q$ are in $C(\Q_p)$ and $\omega$ is a one-form of the second kind. We give a brief outline here, and for more details refer the reader to \cite{BBK}.

To evaluate Coleman integrals, we work locally in each residue disc, where the Coleman integral can be expressed as a convergent power series in a parameter, called a \emph{local coordinate}. If $K(C)$ is the function field of $C$, then our local coordinate is simply a uniformizing parameter $t$ for $K(C)_P$, that is, an element $t \in K(C)_P$ such that the valuation of $t$ is one.

Our residue discs come in two flavors: those that are the residue disc of a Weierstrass point, the \emph{Weierstrass residue discs}, and those that are not, the \emph{non-Weierstrass} residue discs. The \emph{infinity disc} which contains the point at infinity is a Weierstrass residue disc in the odd degree hyperelliptic case.

For a hyperelliptic curve $C/K$ we have differentials $\omega_i = (x^i / 2y) \d{x}$, for $i = 0, \dots, 2g-1$ such that their classes $[\omega_i]$ generate $H^1_{dR}(C)$. For any one-form $\omega$, we can write $\omega = \sum_{i = 0}^{2g-1} a_i \omega_i + \d f$ for some $a_i \in K$ and exact differential $\d f$, and then use linearity to break up an integral.

Assume that $P$ and $Q$ are in the same residue disc, and that this disc does not contain a pole of $\omega$.
To compute a Coleman integral between {$P$ and $Q$, called a \emph{tiny integral}, we simply write $P,Q,$ and $\omega$ in a local coordinate $t$, $P(t)$, $Q(t)$, $\omega(t)$, and integrate formally:
\[ \int_P^Q \omega = \sum_{i = 0}^{2g-1}a_i \int_P^Q  \omega_i + \int_P^Q \d f = \sum_{i=0}^{2g-1} a_i\int_{P(t)}^{Q(t)} \frac{x(t)^i}{2y(t)} \frac{d x(t)}{d t} \d t + f(Q) - f(P).\]
To compute a Coleman integral between $P$ and $Q$ in different non-Weierstrass residue discs we fix Teichm\"uller points $P\p$ in the residue disc of $P$ and $Q\p$ in the residue disc of $Q$ such that $\phi(P\p) = P\p$ and $\phi(Q\p) = Q\p$. By the previous algorithm we can compute tiny integrals $\int_P^{P\p} \omega$ and $\int_Q^{Q\p} \omega$ between our desired points and the Teichm\"uller points. Applying the change of variables property of the integral, $\int_{P\p}^{Q\p} \phi^*\omega = \int_{\phi(P\p)}^{\phi(Q\p)} \omega $. Thus we compute $\phi^* \omega_i$ for each basis vector $\omega_i$. Then $\phi^*$ acts linearly on cohomology, so we can represent $\phi^*$ as a matrix $M$. However, $\phi^* \omega_i$ is only cohomologous to a linear combination of $\omega_i$, so to represent the action of $\phi^*$ as a matrix, we write \[\phi^* \omega_i =  \sum_{j=0}^{2g-1} M_{ij} \omega_j + df_i.\] This gives us a linear system, and we can solve \[ (M-I) \left[\begin{array}{c} \vdots \\ \int_{P\p}^{Q\p}\omega_i\\  \vdots    \end{array}\right] =  \left[\begin{array}{c} \vdots\\ \sum_{i=0}^{2g-1} f_i(P\p) - f_i(Q\p)\\ \vdots \end{array} \right]\]for the value of $\int_{P\p}^{Q\p} \omega_i .$ By linearity, we can compute $\int_{P\p}^{Q\p} \omega .$
The value of the original integral is simply the sum $\int_P^{P\p} \omega + \int_{P\p}^{Q\p} \omega + \int_{Q\p}^Q \omega$. Finally, for the case of $P$ or $Q$ in a Weierstrass disc, we use the fact that if $P$ is a Weierstrass point $\int_P^Q \omega = 1/2 \int_{\iota(Q)}^Q \omega$ where $\iota$ is the hyperelliptic involution (Lemma 16 of \cite{BBK}).

\section{Algorithm}
In this section we describe our implementation of the Chabauty-Coleman method for our case of interest. For a more general description of this method, see for instance \cite{McCallumPoonen}. Our code is available at \cite{marisachigit}.

We require $C$ to be a hyperelliptic curve of genus $3$  given by an odd degree model
\[ C : y^2 = F(x),\]
where $F(x) \in \Q[x]$ is a monic polynomial of degree $7$. Every such model has a unique rational point at infinity, that we denote $\infty$ and take as our basepoint for the Abel-Jacobi embedding $\iota: C \hookrightarrow J$. Moreover, we  assume that the Jacobian $J$ of $C$ has Mordell-Weil rank $0$ over $\Q$. Let $p \geq 7$ be the smallest prime such that $C$ has good reduction mod $p$ (note that $p$ does not divide the leading coefficient of $F$ since we assume that $F$ is monic).

Our goal is to provably compute the set of rational points of $C$. We first use Magma to compute $C(\Q)_{\text{known}}$ to be the set of all rational points in $C$ of naive height bounded by $10^5$. Having this set reduces the number of Coleman integrations required by the algorithm. 

The input of our algorithm consists of the hyperelliptic curve $C$, the prime $p$ and the set $C(\Q)_{\text{known}}$. The output is a finite subset of $C(\Q_p)$ containing $C(\Q)$, which is returned as three separate subsets:
\begin{itemize}
  \item The set $C(\Q)$ of rational points of $C$.
  \item The set of points $Q$ in $C(\Q_p)\setminus C(\Q)$ such that $[Q-\infty] \in J(\Q_p)$ is $2$-torsion.
  \item The set of points $Q$ in $C(\Q_p)\setminus C(\Q)$ such that $[Q-\infty] \in J(\Q_p)$ is an $n$-torsion point for some $n > 2$.
\end{itemize}

Note that since we are assuming that $J$ has rank $0$, no point in $C$ can give rise to a point of infinite order in $J$. Moreover, we remark that under this rank hypothesis we automatically know that the annihilator of $J(\Q)$ under the integration pairing is spanned by $\omega_i = x^i/2y \d{x}$ for $i = 0,1,2$, which simplifies Step 2 of our algorithm.

We use essentially the same precision analysis as in \cite{CCexp}. Proposition 3.11 in \cite{CCexp} can be easily adapted to show that it is enough to set $N = 2p+4$ to the $p$-adic precision and $M = 2p+1$ to be the $t$-adic precision. We give a basic outline the adaptations needed here (the reader should replace the statements in their paper about the annihilating differentials $\alpha, \beta$ with the annihilating differentials $\omega_0, \omega_1, \omega_2$).

First note that assuming \cite[Proposition~3.11]{CCexp}, $N = 2p+4$ and $M = 2p+1$ is sufficient for $p$ prime $p>2g = 6$ (this follows from Riemann-Roch, see \cite[Remark~3.12]{CCexp}). The proof of this proposition depends on two earlier statements: Lemmas 3.10 and 3.9, which do not depend on the rank of the curve or the exact set of annihilating differentials.

While the lemmas and propositions in \cite[Section~3.1]{CCexp} that bound the number of zeros in each residue disc of the curve remain true with minimal modification to the proofs, we did not use these bounds to rule out residue discs in our implementation.

An outline of our algorithm is given in Algorithm \ref{main-alg}. Below we explain its main steps.

\begin{algorithm}
\caption{Chabauty-Coleman method for a genus $3$ hyperelliptic curve with rank $0$ Jacobian} \label{main-alg}

\begin{algorithmic}[1]
  \Function{Chabauty-Coleman}{$C$, $p$, $C(\Q)_{\text{known}}$}
  \State Set the $p$-adic precision $N$ to $2p+4$.
  \State Set the $t$-adic precision $M$ to $2p+1$.
  \State Initialize found-points := empty list. 
  \For {each $\overline{P} \in \overline{C}(\Ff_p)$ up to the hyperelliptic involution $\iota$}
  \State Compute $f_0, f_1$ and $f_2$ in local coordinates.
  \For {each point $Q \in C(\Q_p)$ corresponding to a common zero of the $f_i$}
  \State Add $Q$ and $\iota(Q)$ to found-points.
  \EndFor
  \EndFor
  \State Classify found-points into three lists: $\Q$-points, non-rational $2$-torsion-points, and non-rational higher-torsion-points.\\
  \Return $\Q$-points, non-rational $2$-torsion-points, and non-rational higher-torsion-points.
  \EndFunction
  \end{algorithmic}
\end{algorithm}

\textit{Step 1 (Required precision.)} We need to choose the $p$-adic precision $N$ and the $t$-adic precision $M$ to guarantee that, in Step 3, we will obtain all the roots of $f_i(pt)$ in $\Z_p$. Set $N = 2p+4$ and $M = 2p+1$.

\textit{Step 2 (Annihilator.)} A basis of the space of differentials $H^0(C_{\Q_p}, \Omega^1)$ is given by $\{\omega_0, \omega_1, \omega_2 \}$, where $\omega_i = (x^i/2y)\d{x}$.
For each $i = 0, 1, 2$, define
$$f_i(z) = \int_\infty^z \omega_i,$$
where $\infty$ denotes the point at infinity. The functions $f_i(z)$ are zero on all rational points of $C$, but not identically zero.

\textit{Step 3 (Searching in residue discs.)} For each point $\overline{P}$ in $\overline{C}(\Ff_p)$, we compute the set of $\Q_p$-rational points $P$ reducing to $\overline{P}$ such that $f_i(P) = 0$ for $i=0, 1, 2$. To perform this computation, we consider two different cases:

\begin{enumerate}[(i)]
\item If there is a point $P \in C(\Q)_\text{known}$ reducing to $\overline{P}$, then we expand each $\omega_i$ in terms of a uniformizer $t$ at $P$ and we formally integrate to obtain three power series $f_i(t)$, that parametrize the integrals of the $\omega_i$ between $P$ and any other point in the residue disc.
\item Otherwise, we start by finding a $\Q_p$-point $P$ reducing to $\overline{P} = (\overline{x_0}, \overline{y_0})$ (note that $\overline{P}$ cannot be $\overline{\infty}$ in this case). If $\overline{y_0}=0$ we can take $P = (x_0, 0)$ where $x_0$ is the Hensel lift of $\overline{x_0}$ to a root of $F(x)$.
Otherwise,  we can take $P = (x_0 , y_0)$ where $x_0$ is any lift of $\overline{x_0}$ to $\Z_p$ and $y_0$ is obtained from $\overline{y_0}$ by applying Hensel’s Lemma to $y^2 = F(x_0 )$.
Then we set $f_i(t) = \tilde{f}_i(t) + \int_\infty^P \omega_i,$ where each $\tilde{f}_i(t)$ parametrizes the integral of $\omega_i$ between $P$ and any other point in the residue disc. 
\end{enumerate}

To provably compute the set of common zeros of $f_0, f_1$ and $f_2$ to a desired precision, we need to impose the condition that at least one of the $f_i$ has only simple roots. We check this requirement by computing the discriminant of the truncations of the power series to $t$-adic precision $M$. If all the $f_i$ have multiple roots, we run the algorithm using a bigger prime $p$ for which they do not.

Say $f_{i_0}$ has no double roots. Then we use the PARI/GP function \textsf{polrootspadic} to compute its roots, truncating it first to $t$-adic precision $M=2p+1$ and computing the coefficients up to $p$-adic precision $N=2p+4$. For each root that lies in $p\Z_p$, we check whether it is also a root of the other two $f_i$; if so, it corresponds to a point $P \in C(\Q_p)$ lying over $\overline{P}$.

\textit{Step 4 (Identifying the rational points.) }Now, for each of the points $Q$ found in Step 3, we attempt to reconstruct $Q$ as a $\Q$-rational point, using Sage. If this is not possible, then $[Q - \infty]$ must be a torsion point in $J(\Q_p)$, because $J$ has rank $0$. If $Q$ is a Weierstrass point, it will give rise to a $2$-torsion point in the Jacobian; otherwise we classify it as a higher order torsion point.

\textit{Step 5 (Identifying higher order torsion points.)} For each of the points $Q$ that is not identified as a rational or Weierstrass point in Step 4, we use Sage's \textsf{algebraic\_dependency} function to attempt to find the minimal polynomial of the $x$-coordinate and reconstruct the point as a point with coordinates defined over an extension of $\Q$. We then use Magma to compute the torsion order of $[Q - \infty]$.

\section{Analysis and Examples}

We ran the Sage implementation of our algorithm on a database of $5870$ hyperelliptic curves of genus $3$ with Jacobian of rank $0$ taken from the list \cite{g3rk0}. This list \cite{g3rk0} is a subset of Sutherland's database of genus 3 hyperelliptic curves of discriminant bounded by $10^7$ \cite{g3database}. It was created by first filtering for curves with odd-degree models and then filtering for curves that are provably rank 0. The rank computation was done with the Magma function \texttt{RankBounds}, which returns an upper and lower bound for the rank of the curve through descent computations. For curves in \cite{g3rk0}, the upper bound and the lower bound returned by \texttt{RankBounds} agreed and equaled 0. We computed a set of monic models for each curve in \cite{g3rk0} for the purpose of our computations.

In $23$ of these curves, there was at least one $\Q_p$-point for which the three power series $f_0(z), f_1(z) \text{ and } f_2(z)$ had double roots at the first prime $p \geq 7$ of good reduction. To analyze these curves, we replaced $p$ by a bigger prime for which at least one of the power series had only simple roots.

Our implementation proves that, for each of the studied curves, the entire set of rational points equals the set of rational points of naive height less than $10^5$. In fact, we proved that, with respect to the odd-degree monic models used in the computation, the maximum among the global heights of rational points in those models is $30.7611440827071$, reached only by points at three curves. Or, if we translate this to the odd degree models in \cite{g3rk0}, those three curves have points of global height $4.39444915467244$, while the heights of the rational points in the rest of the curves is bounded above by $3$. \footnote{The global height is the absolute logarithmic height of the point, which is the maximum of the absolute logarithmic heights of its coordinates. For a rational point $\frac{n}{d}$, this height is computed as $\max(\log(|n|),\log(|d|)).$}

Next, we show in
Figure \ref{fig:points} how many of the curves in our database have a certain number of rational points. We observe that the maximum number of points is six, and that a vast majority of the curves have three or fewer rational points.

\begin{figure}[h!]
  \begin{center}
    \includegraphics[width=0.5\textwidth]{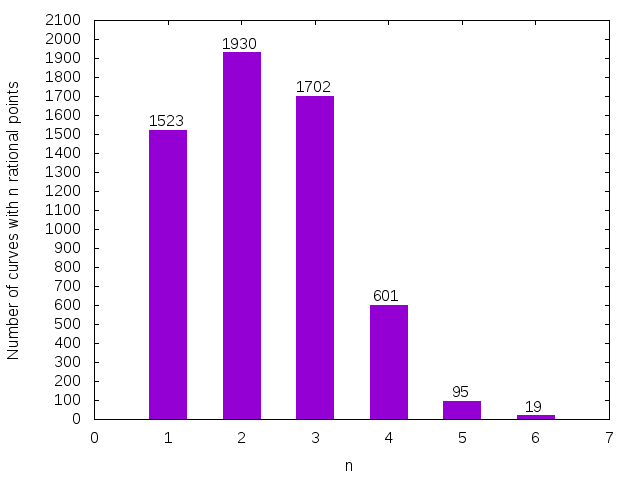}
    \caption{Number of curves in the database with $n$ rational points}
    \label{fig:points}
  \end{center}
\end{figure}
%\vspace{-0.5cm}

As mentioned on the introduction, the Chabauty--Coleman sets of these curves agreed with the set of rational points on $2783$ cases, contained only rational points and Weierstrass points on $3070$ cases, and contained points on the curve giving rise to higher order torsion points on the Jacobian on $17$ cases, $13$ of which also presented non-rational Weierstrass points . In each of these $17$ curves, there was exactly one pair of points corresponding to higher order torsion points, of order $8$ in ten cases, order $12$ in six cases, and order $18$ in one case, which we present in Example $2$. We made these $17$ higher torsion examples available on \cite{marisachigit}.

We conclude this section by showing how the algorithm works on a particular curve, giving an example where the algorithm detects a torsion point defined over a number field, and exhibiting an example curve where the Stoll bound is sharp.

\begin{eg}
          
Consider the hyperelliptic curve 
\begin{align*}
C: y^2 = x^7 &- 37024x^6 + 3134464x^5 - 101220352x^4 + 1613758464x^3 -\\ &-13656653824x^2 + 59055800320x - 103079215104,
\end{align*}
with minimal discriminant $2168084$ and conductor $1084042$.

Using Magma, we find that the set of rational points of $C$ with height bounded by $10^5$ is
\[C(\Q)_{\text{known}} = \{\infty, (32, 0)\}.\]
Since $C$ has good reduction modulo $7$, we run our code using this prime. The points of $C(\mathbb{F}_7)$ are as follows:
\[ \{ \overline{\infty} , \overline{ (0 , \pm 4 ) } , \overline{ (1 , \pm 5 )}, \overline{ (2 , \pm 6 )}, \overline{ (4 ,0) }, \overline{ (6 ,\pm 2 )} \}. \]

After Hensel lifting each of these points to a point of $C(\Q_7)$, we write $f_0(z)$, $f_1(z)$, and $f_2(z)$ in local coordinates and find their common zeros. This yields the following results:
\begin{center}
  \renewcommand{\arraystretch}{1.5}
  \begin{tabular}{ |c|c| } 
    \hline
    disc & common roots of $f_0(z)$,  $f_1(z)$ and  $f_2(z)$ \\[3pt] \hline
    $\overline{\infty}$ & $\infty$ \\
    \hline
    $\overline{ (0 , \pm 4 ) }$  & no common roots \\[2pt]
    \hline
    $\overline{(1,\pm 5)} $& no common roots \\[2pt]
    \hline
    $\overline{(2 ,\pm 6)}$ & no common roots \\ [2pt]
    \hline
    $\overline{ (4 ,0) }$ &  $(32, 0)$ \\[2pt]
    \hline
    $\overline{ (6 ,\pm 2 )} $ & no common roots \\ [2pt]
    \hline
  \end{tabular}
  \renewcommand{\arraystretch}{1}
\end{center}
Therefore, we have shown that $C(\Q)= \{ \infty, (32,0)\}$.
\end{eg}

\begin{eg}
Let \[C: y^2 = x^7 + \frac{5}{64}x^6 - \frac{51}{256}x^5 - \frac{243}{4096}x^4 - \frac{53}{8192}x^3 - \frac{27}{65536}x^2 - \frac{1}{65536}x - \frac{1}{4194304}\] be the hyperelliptic curve with minimal discriminant 6856704 and conductor 6856704. Using the algorithm described in the paper, we can determine that this curve has only one rational point, the point at infinity. However, we also detect a pair of $7$-adic points with $x$-coordinate

\begin{center}
 \begin{tabular}{c} 
$6 + 6\cdot7^2 + 6\cdot7^4 + 6\cdot7^6 + 6\cdot7^8 + 6\cdot7^{10} + 6\cdot7^{12} + 6\cdot7^{14} + 6\cdot7^{16} + O(7^{18})$.\\

\end{tabular}
\end{center}

Using sage's \textsf{algebraic\_dependency} function, we identify this pair of points as $(-1/8,\pm \sqrt{-3}/2^{11})$. The points $[(-1/8,\pm \sqrt{-3}/2^{11}) - \infty]$ have order 18 on the Jacobian over $\Q(\sqrt{-3})$. The Chabauty-Coleman algorithm detects these torsion points at the prime $p =7$ because $7$ splits in this field.

\end{eg}

\begin{eg}
When $p>2g$ is a prime of good reduction and $r<g-1$, then Stoll \cite[Corollary~6.7]{StollBound} has improved the bound given by Coleman to give $\# X(\Q) \leq \#X(\Ff_p) + 2r$. In our case, this shows that $\#X(\Q) \leq \#X(\Ff_p)$. In our database of 5870 curves, 94 achieve equality for Stoll's bound.

For example, the curve
% \[y^2 = 8x^7 - 16x^5 - 7x^4 + 4x^3 + 6x^2 + 4x + 1,\]
\[ y^2 =  x^7 - \frac{1}{2^5}x^5 - \frac{7}{2^{12}}x^4 + \frac{1}{2^{13}}x^3 + \frac{3}{2^{17}}x^2 + \frac{1}{2^{19}}x + \frac{1}{2^{24}},\]
 which has minimal discriminant 4089600 and conductor 170400, has $6$ rational points of height less than $10^5$: 
 %\[\{\infty, (0 , 1 ), (0 , -1 ), (1 , 0 ),(-1 , 0 ) , (-1/2 , 0 )\}\]
 \[\left\{ \infty, \left(0, \frac{1}{4096}\right), \left(0, -\frac{1}{4096}\right), \left(\frac{1}{8}, 0\right), \left(-\frac{1}{8}, 0\right), \left(-\frac{1}{16}, 0\right)  \right \}\]
  as well as $6$ points over $\Ff_7$:
\[\{\infty, (0 , 1), (0 , 6 ),(1 , 0 ), (6 ,0 ),(3 , 0)\}.\]
By Stoll's theorem, we know that there cannot be any more rational points.

\end{eg}

\end{document}